\documentclass{amsart}
\usepackage[utf8]{inputenc}

\usepackage[asymmetric,vmargin=3cm,hmargin=3cm,marginpar=2cm]{geometry}

\usepackage{amsmath,amssymb,amsfonts,amsthm,amscd,graphicx,psfrag,epsfig,fge}
\usepackage{bbm}
\usepackage{shuffle}

\usepackage{enumitem}
\setlist[enumerate]{itemsep=1ex,leftmargin=2em, label={\rm(\roman*)}}
\setlist[itemize]{itemsep=1ex,leftmargin=2em}

\usepackage{xcolor}
\usepackage{hyperref}
\usepackage{mathtools}

\usepackage{thmtools}

\declaretheorem[
style=plain,
name=Theorem,
numbered=yes,
refname={Theorem,Theorems},
Refname={Theorem,Theorems}
]{theorem}

\declaretheorem[
style=plain,
name=Lemma,
numberlike=theorem,
refname={Lemma,Lemmas},
Refname={Lemma,Lemmas }
]{lemma}

\allowdisplaybreaks

\title[A property of trivalent graphs related to equidissections
]{A property of trivalent graphs related to equidissections}

\author[D. Rudenko]{Daniil Rudenko}

\date{\today}

\begin{document}


\begin{abstract}
Monsky proved that a square cannot be dissected into an odd number of triangles of equal area. Stein conjectured that the same holds for any polygon whose edges can be paired into parallel and equal-length segments. We prove Stein’s conjecture under an assumption that all triangle vertices have rational coordinates. Our result is derived from a more general property of trivalent graphs equipped with a $\mathbb{Q}^2$-valued flow.
\end{abstract}

\maketitle

\section{Introduction}
Consider a finite trivalent graph $\Gamma$ without loops or multiple edges. Let $V$ and $E$ denote its sets of vertices and edges, respectively. We denote by $\vec{E}$ the set of \emph{oriented edges} of $\Gamma$—that is, pairs consisting of an edge and its orientation. Let $A$ be an abelian group. An \emph{$A$-flow} on $\Gamma$ is a function $f\colon \vec{E} \to A$ such that:
\begin{enumerate}[label=(\roman*)]
  \item $f(e) + f(e') = 0$ for every pair of opposite oriented edges $e, e' \in \vec{E}$,
  \item $f(e_1) + f(e_2) + f(e_3) = 0$ for every triple of oriented edges $e_1, e_2, e_3 \in \vec{E}$ with common source vertex.
\end{enumerate}

Now consider a trivalent graph $\Gamma$ with a flow $f$ valued in $A = \mathbb{Q}^2$. Let $\omega \in (\Lambda^2 \mathbb{Q}^2)^\vee$ be the standard area form given by
$\omega\bigl((x_1, y_1), (x_2, y_2)\bigr) = x_1 y_2 - x_2 y_1.$
The \emph{weight} of a vertex $v \in V$ is defined as
\[
m(v) = \nu_2\Bigl(\omega\bigl(f(e_1), f(e_2)\bigr)\Bigr),
\]
where $\nu_2 \colon \mathbb{Q} \to \mathbb{Z} \cup \{\infty\}$ denotes the $2$-adic valuation, and $e_1, e_2$ are any two of the three oriented edges with the source $v$. This quantity is well-defined due to the antisymmetry of $\omega$ and the flow condition. The \emph{weight of the graph} $\Gamma$ with the flow $f$ is the minimum of the vertex weights:
\[
m(\Gamma,f) = \min_{v \in V} m(v).
\]

Our main result is the following.

\begin{theorem} \label{MainTheorem} 
Let $\Gamma$ be a trivalent graph with a flow $f\colon \vec{E} \to \mathbb{Q}^2$. Then the number of vertices $v \in V$ such that $m(v) = m(\Gamma)$ is even.
\end{theorem}

Theorem \ref{MainTheorem} is a result of our attempts to understand the proof of Monsky's theorem (\cite{Mon70}), which states that a square cannot be dissected into an odd number of triangles of equal area. Here a dissection of a simply connected polygon $P$ is a finite set of triangles with disjoint interiors that cover $P$. Monsky's theorem is known for its elegant proof which uses a $2$-adic valuation and Sperner's lemma. Later Monsky proved a much more difficult result: the same statement holds for any centrally-symmetric polygon (\cite{Mon90}).

Stein conjectured (\cite{Ste00}) that a similar statement holds for a broader class of polygons, which he called  special. Consider a simply-connected polygon $P$ and fix an orientation of its boundary; each side of $P$ determines a vector in $\mathbb{R}^2$. The polygon $P$ is said to be {\it special} if, for every line $l$, the sum of the vectors corresponding to the sides of $P$ that are parallel to $l$ is zero. Stein conjectured that a special polygon cannot be dissected into
an odd number of triangles of equal area. We prove that conjecture under an additional rationality assumption:

\begin{theorem}[Stein's conjecture for rational dissections] \label{MainCorollary} A special polygon $P\subseteq \mathbb{R}^2$ cannot be dissected into an odd number of triangles of equal area
in such a way that all coordinates of the vertices of the triangles in the dissection are rational numbers.
\end{theorem}

We derive Theorem \ref{MainCorollary} from Theorem \ref{MainTheorem}. For a given dissection of a polygon, we consider a trivalent graph $\Gamma$, which is obtained from the dual graph of the dissection by adding extra degenerate triangles and connecting vertices corresponding to triangles adjacent to the boundary of $P$. The flow assigns to an edge of $\Gamma$ the side vector of the corresponding triangle of the dissection. If the triangles in the dissection have the same area $A$, the graph $\Gamma$ has vertices of weight $\nu_2(A)\in \mathbb{Z}$ and $\nu_2(0)=\infty$. Theorem \ref{MainTheorem}  implies that the number of vertices of weight $\nu_2(A)$ must be even.

\parskip 9pt
{\bf Acknowledgments}.
I am grateful to Sergei Tabachnikov for introducing me to the topic of equidissections, and to Nikolai Mnev, whose guidance and support were essential to the completion of this article. I also thank Dmitri Krachun for many helpful discussions.

\section{Proof of Theorem \ref{MainTheorem}}

We begin with some preliminaries. Let $\mathbb{Z}_{(2)}$ be the localization of $\mathbb{Z}$ at the prime ideal $(2)$; elements of $\mathbb{Z}_{(2)}$ are rational numbers with an odd denominator. This is a discrete valuation ring with the discrete valuation $\nu_2$, the maximal ideal $(2)$, and the residue field $\mathbb{F}_2$. A lattice over $\mathbb{Z}_{(2)}$ is a free module of finite rank. Consider a lattice $L$ and its sublattice $L'$. The index $(L:L')$ is finite if and only if $L$ and $L'$ have the same rank; it equals $2^{\ell(L,L')}$ where $\ell$ is the length. If we choose a basis $e_1,\dots, e_n$ of $L$ and $e_1',\dots, e_n'$ of $L'$ such that $e_i'=\sum a_{ij} e_j$, then the length $\ell(L,L')$ is equal to the $2$-adic valuation of the determinant of the matrix $(a_{ij})$. 

Let $\Gamma$ be a trivalent graph with a $\mathbb{Q}^2$-valued flow. For every vertex $v\in \Gamma$, we denote by $L(v)$ a $\mathbb{Z}_{(2)}$-submodule of $\mathbb{Q}^2$ spanned by vectors $f(e_1)$ and $f(e_2)$, for any two edges $e_1$ and $e_2$ with the source~$v$. The module $L(v)$ is free of rank at most two. We say that a flow is {\it integral} if $f(e)\in \mathbb{Z}_{(2)}^{2}$ for every oriented edge $e\in \vec{E}$. If the flow is integral, the lattice $L(v)$ is contained in $\mathbb{Z}_{(2)}^{2}$ for every vertex $v\in V$ and 
$m(v)=l\bigl(\mathbb{Z}_{(2)}^{2},L(v)\bigr).$

 Consider the projection 
$p\colon \mathbb{Z}_{(2)}^2\longrightarrow \mathbb{F}_2^2.$
A vectors $u\in \mathbb{Z}_{(2)}^2$ is called {\it even} if $p(u)=0$, and {\it odd} otherwise.
If $f$ is an integral flow, we call an edge $e$  {\it even} if $f(e)$ is even, and {\it odd} if $f(e)$ is  odd. In this case, the weight of a vertex $v$ is equal to zero if and only if we can order the edges with source $v$ so that $p(f(e_1))=(0,1)$,  $p(f(e_2))=(1,0)$, and $p(f(e_3))=(1,1)$. If the weight of  $v$ is strictly positive, then at least one of the edges $e_1, e_2, e_3$ is even. Since $f$ is a flow, it follows that, in this case, either one or three of these edges is even. Thus, the number of vertices of weight zero has the same parity as the number of edges of $\Gamma$. Since $\Gamma$ is trivalent, this number is even. We have proven the following lemma:
\begin{lemma} \label{LemmaMain} Let $\Gamma$ be a trivalent graph with an integral flow $f$. The number of vertices of $\Gamma$ of weight zero is even.
\end{lemma}

Now we are ready to prove Theorem \ref{MainTheorem}. Its statement  is invariant under rescaling, so we may assume that the flow $f$ is integral. We argue by contradiction and choose a trivalent graph $\Gamma$ with an integral flow which does not satisfy Theorem \ref{MainTheorem} with minimal weight among such counterexamples. If the weight $m(\Gamma,f)$ is equal to zero, we get a contradiction with Lemma \ref{LemmaMain}. Therefore, we may assume that the weight $m(\Gamma)$ is strictly positive. 

In this case, every vertex is the source of an even number of odd edges, so the odd edges form a union of disjoint cycles; we will call these cycles {\it odd}. The graph $\Gamma$ must contain at least one odd cycle, because otherwise all edges of $\Gamma$ would be even, and the same graph with the flow $\frac{1}{2}f$ would be a counterexample of smaller weight. Consider an odd cycle $v_1, v_2,\dots, v_k, v_{k+1}=v_1$ and denote the oriented odd edge from $v_i$ to $v_{i+1}$ by $e_{i,i+1}$.  Let $e_i$ be the even edge with source $v_i$, and denote its target by $v_i'$. Note that some of the vertices $v_i'$ may coincide.  Clearly, $f(e_1)+\dots+f(e_k)=0$. 
Finally, denote the lattice $L(v_i)$ by $L_i$. 

\begin{lemma} Let $v_1, v_2,\dots, v_k, v_{k+1}=v_1$ be an odd cycle in a trivalent graph with an integral flow. There exists $i\in \{1,\dots, k\}$ such that $L_j\subseteq L_i$ for any $j\in \{1,\dots, k\}$.
\end{lemma}
\begin{proof}
 We start with an observation. Let $u\in \mathbb{Z}_{(2)}^2$ be an odd vector. Then the lattices contained in  $\mathbb{Z}_{(2)}^2$ and containing $u$ are linearly ordered by inclusion, because they are in one-to-one correspondence with the submodules of $\mathbb{Z}_{(2)}^2/\langle u\rangle \cong \mathbb{Z}_{(2)}$. Now, choose a lattice $L_i$ which is not strictly contained in any lattice $L_j$; we may assume that $i=1$. Assume that there exists a lattice $L_j$ not contained in $L_1$, and choose the smallest such $j$. Then $L_{j-1}\subseteq L_1$. Since  $L_{j-1}$ and $L_j$ contain the vector $f(e_{j-1,j})$ with an odd coordinate, we have $L_j\subseteq L_{j-1}$ or $L_{j-1}\subseteq L_j$. In the first case, 
$L_j\subseteq L_{j-1} \subseteq L_1$,
which contradicts our assumption.  In the second case,  lattices $L_1$ and $L_j$ contain the lattice $L_{j-1}$, and, consequently, contain an odd vector. Since $L_1$ is not strictly contained in $L_j$, we must have $L_j \subseteq L_1$, which contradicts our assumption.
\end{proof}

\begin{lemma} \label{Lemma3} Let $v_1, v_2,\dots, v_k, v_{k+1}=v_1$ be an odd cycle in a trivalent graph with an integral flow and let $L$ be the maximal lattice among $L_1,\dots, L_k$. Then the number of indices $i\in \{1, \dots, k\}$ for which $L_i=L$ is even.
\end{lemma} 
\begin{proof}
    Consider a subgraph of $\Gamma$ induced by the vertices $v_1,\dots, v_k$ and $v_1',\dots, v_k'$. This subgraph can be extended to a trivalent graph $C$ by adding an edge $e_{i,i+1}'$ from $v_i'$ to $v_{i+1}'$. There exists a unique flow $f_C$ on $C$  which coincides with $f$ on edges $e_{i,i+1}$ and $e_i$ and vanishes on the edge $e_{1,2}'$. Let $L^{0}$ be the lattice spanned by the even vectors in $L$. From the fact that $L$ is a proper sublattice of $\mathbb{Z}_{(2)}^2$ it is easy to see that the index of $L^{0}$ in  $L$ equals to two. We have $e_{i,i+1}'\in L^{0}$ for $i\in \{1,\dots,k\}$. So, the lattices at the vertices of $C$ are contained in $L$, while those corresponding to the vertices $v_i'$ are proper sublattices of $L$.
    Consider a linear map $A\in \mathrm{GL}_2(\mathbb{Q})$ sending the lattice  $L$ to  $ \mathbb{Z}_{(2)}^2$. The statement of the lemma follows from Lemma \ref{LemmaMain} applied to the integral flow $f'_C:=A\circ f_C$ on $C$.
\end{proof}

To finish the proof of Theorem~\ref{MainTheorem}, we construct a new flow $f'$ on $\Gamma$ for which all edges are even, $m(\Gamma, f) = m(\Gamma, f')$, and the number of vertices of weight $m(\Gamma, f')$ is odd. Then the graph $\Gamma$ with the flow $\tfrac{1}{2}f'$ would be a counterexample to Theorem~\ref{MainTheorem} with an integral flow and strictly smaller weight than $\Gamma$. This contradicts our choice of $\Gamma$.

We construct the flow $f'$ as follows. If $e$ is an even oriented edge, we set $f'(e):=f(e)$. Every odd edge of $\Gamma$ is a part of some odd cycle. For each odd cycle $v_1,\dots, v_k, v_{k+1}=v_1$, and $i\in \{1,\dots, k\}$, we set $f'(e_{i,i+1})=f(e_{i,i+1})-f(e_{1,2})\in L^0$. Observe that $f'$ is even on every edge of $\Gamma$.  Since vectors $f'(e_{i,i+1})$ are contained in $L^0$, which is a proper sublattice of $L$, we have not created vertices of weight greater than $m$. Finally, Lemma \ref{Lemma3} implies that the number of vertices of weight $m$ in the graph $\Gamma$ with the flow $f$ and with the flow  $f'$ have the same parity. As we have  explained above, the existence of a flow with these properties contradicts our choice of $\Gamma$ and $f$. This completes the proof of Theorem \ref{MainTheorem}.

\section{Proof of Theorem \ref{MainCorollary}}

A {\it degenerate triangle} is a triple of distinct points in $\mathbb{R}^2$ that lie on a common line. A {\it generalized triangle} is either a triangle or a degenerate triangle. An {\it orientation} of a generalized triangle is an ordering of its vertices; this induces an orientation on each of its edges.

Suppose that $P$ is a special polygon dissected into triangles of equal area $A$, and let $\mathcal{P}$ denote the set of all vertices of triangles in the dissection. By adjoining degenerate triangles with vertices in~$\mathcal{P}$, we may assume that any two generalized triangles intersect in either an empty set, a single vertex, or a common edge. Next, by adding points with rational coordinates to the boundary of $P$, we may assume that the sides of the triangles lying on the boundary of $P$ can be paired so that the corresponding oriented vectors in each pair sum to zero. We fix such a pairing and refer to the paired boundary sides as {\it associated}.

We now orient the generalized triangles of the dissection so that:
\begin{enumerate}
    \item if two generalized triangles share a common edge, they induce opposite orientations on that edge;
    \item if an edge lies on the boundary of $P$, then the orientation it inherits from the unique triangle containing it agrees with the orientation of the boundary of~$P$.
\end{enumerate}

Let $V$ be the set of generalized triangles constructed above. Define a graph $\Gamma$, where $\vec{E} \subset V \times V$ consists of pairs of generalized triangles $(T, T')$ that either share a common edge or contain associated boundary sides. Define a function $f \colon \vec{E}\longrightarrow \mathbb{Q}^2$ by assigning to each pair $(T, T')$ the vector corresponding to the side of $T$, with orientation induced by that of $T$. It is easy to see that $\Gamma$ is a trivalent graph and that $f$ defines a $\mathbb{Q}^2$-valued flow on it.

Now consider a vertex of $\Gamma$ corresponding to a non-degenerate triangle $T$, and let $e_1, e_2$ be any two of the three edges with source $T$. The area of $T$ is equal to $\pm\frac{1}{2}\omega(v_1, v_2)$, where $v_1$ and $v_2$ are the vectors assigned to $e_1$ and $e_2$, respectively. It follows that the weight of such a vertex is equal to $\nu_2(A) - 1$. Thus, all vertices corresponding to non-degenerate triangles have the same finite weight. On the other hand, all vertices corresponding to degenerate triangles have infinite weight. By Theorem~\ref{MainTheorem}, it follows that the number of non-degenerate triangles in the dissection is even. This completes the proof of Theorem~\ref{MainCorollary}.

\bibliographystyle{alpha}      
\bibliography{bibliography} 
\end{document}